\documentclass[a4paper,12pt]{article}

\usepackage[frenchb,english]{babel}
\usepackage{amsmath,amssymb,amscd,amsfonts,amsthm}
\usepackage[T1]{fontenc}
\usepackage[latin1]{inputenc}
\usepackage{indentfirst}

\newcommand{\sph}{\mathbb{S}}
\newcommand{\tor}{\mathbb{T}}
\newcommand{\ha}{\mathbb{H}}
\newcommand{\R}{\mathbb{R}}
\newcommand{\finthm}{\hfill $\square$ \vspace{4mm}}

\newtheorem{theorem}{Theorem}
\newtheorem{cor}[theorem]{Corollary}
\newtheorem{lemma}[theorem]{Lemma}
\newtheorem{prop}[theorem]{Proposition}
\theoremstyle{definition}
\newtheorem{definition}[theorem]{Definition}

\newtheorem{example}[theorem]{Example}

\begin{document} 
\sloppy

\author{Mathias Zessin}
\title{On contact tops and integrable tops}
\date{}
\maketitle \vspace{-1cm}

{\small \begin{center}
        Mathematisches Institut der Universität zu Köln\\
        Weyertal 86-90 \\
        D-50931 Köln \\
        \underline{e-mail:} mzessin@math.uni-koeln.de
        \end{center}

Keywords~: Contact geometry, Riemannian geometry, contact circles.

MSC 2000~: 57M50, 53D10, 53C21.} \vspace{1mm}

\begin{abstract} In this paper, we introduce a geometric structure called top, which is a trivialized bundle of
plane pencils over a Riemannian 3-manifold, defined as the set of kernels of a circle of 1-forms (e.g. of contact
and integrable forms) with particular properties with respect to the metric. We classify the manifolds which admit
tops and we describe the associated metrics.
\end{abstract}

\section{Definitions and main results} \label{introduction}
The central object of this paper are tops\footnote[1]{Here a top is not the integrable system from Hamiltonian
mechanics, as it is in the usual sense. We use the same terminology, because this word describes so well the
structure we introduce.}. These new structures are defined as trivialized plane pencil bundles over a
3-dimensional Riemannian manifold with some regularity properties with respect to the metric. More precisely, we
define them as the set of kernels of a circle of 1-forms.

There are two different motivations for this work. First, trivialized plane pencil bundles appear as the
geometric structures associated to contact circles, that is, to circles of contact forms, in the same way as
contact structures are associated to contact forms. We present these structures in the context of Riemannian
manifolds and characterize certain families of contact circles which have special metric properties.

The second aim is to answer the following question~: In which extent can the properties of a particular example
on $\sph^3$ (see Examples \ref{ex sphere} and \ref{ex contact tops}) with the usual metric be found on other
manifolds and in relation with other metrics? \vspace{4mm}

Under some additional coorientability assumptions, we associate a circle of differential 1-forms to a trivialized
plane pencil bundle on a 3-manifold $M$, that is, a family
$\mathcal{S}^1\{\omega_1,\,\omega_2\}=\{\cos\theta\,\omega_1+\sin\theta\,\omega_2,\theta\in\left[0,2\pi\right[\}$
generated by two 1-forms $\omega_1$ and $\omega_2$ which are linearly independent in every point of $M$. Given a
two-dimensional subbundle $\Omega$ of the space of differential 1-forms on $M$, the kernels of the non-zero
elements of $\Omega$ define a plane pencil in the tangent space at any point $x \in M$. Moreover, all multiples
of a given form by a positive function have the same kernel. So it is enough to consider the unit circle of
$\Omega$ to realize all the planes of the plane pencil as kernels of elements of this family of 1-forms.

For example, a trivialized pencil bundle can be defined by a contact circle (see section \ref{def contact} for a
precise definition). In that case, the pencil bundle is given by the set of the contact structures associated to
the elements of the contact circle. \vspace{4mm}

All objects in this paper are supposed to be smooth. \vspace{4mm}

Let $\mathcal{F}$ be a plane pencil bundle obtained as the set of kernels of a circle of 1-forms
$\mathcal{S}^1\{\omega_1,\,\omega_2\}=\{\omega^\theta:=\cos\theta\,\omega_1+\sin\theta\,\omega_2,
\theta\in\left[0,2\pi\right[\}$. Then the map
\begin{eqnarray*}
\begin{array}{rccl}
\tau: & \mathcal{F} & \longrightarrow & M\times \mathbb{P}^1(\mathbb{R}) \\
& \ker(\omega_x^\theta) & \longmapsto & (x,\theta\hspace{-3mm}\mod\pi)
\end{array}
\end{eqnarray*}
is a trivialization of $\mathcal{F}$. The plane bundles $\xi_\theta=\ker(\omega^\theta)$ are obtained by the
"inverse map"
\begin{eqnarray*}
\mathbb{P}^1(\mathbb{R}) & \longrightarrow & \mathcal{F} \\
\theta & \longmapsto & \xi_\theta
\end{eqnarray*}
such that $\pi_2(\tau(\xi_\theta))=\theta$. This family is double-covered by the $\sph^1$-family of plane bundles
obtained as the kernels of the defining circle of 1-forms.

\begin{definition} --- \label{ind pencil bundle}
A trivialized plane pencil bundle which is defined by the set of kernels of the elements of a circle of 1-forms
$\mathcal{S}^1\{\omega_1,\,\omega_2\}$ will be called an \emph{indexed pencil bundle}.

An indexed pencil bundle is called \emph{contact pencil bundle} if for all $\theta\in \mathbb{P}^1(\mathbb{R})$,
the plane bundle $\xi_\theta$ is a contact structure, and \emph{integrable pencil bundle} if all these plane
bundles are integrable.

The line bundle $L_x=\cap_{\theta\in\mathbb{P}^1(\mathbb{R})} (\xi_\theta)_x$ is called the \emph{axis bundle} of
the indexed pencil bundle.
\end{definition}

\noindent {\bf Remark} --- If $\mathcal{F}$ is an indexed pencil bundle, then every plane bundle which is
obtained as the kernel of some element of the defining circle of 1-forms is coorientable. On the other hand, a
trivialized plane pencil bundle with the property that all plane bundles $\xi_\theta$ with
$\pi_2(\tau(\xi_\theta))=\theta$ are coorientable can be obtained via a circle of 1-forms, given by
$\omega^\theta=g(X_\theta,\cdot)$, where $X_\theta$ is a differentiable family of vector fields such that
$(X_\theta)_x \bot (\xi_\theta)_x, \forall x\in M, \forall \theta \in \left[0,\pi\right]$, and where $g$ is a
metric on $M$. \vspace{4mm}

A trivialized plane pencil bundle on a Riemannian manifold can have a certain number of geometric regularity
properties. In particular, we can consider the way the plane bundles $\xi_\theta$ rotate about the axis bundle
with respect to the parallel transport along some particular curves. To make this more precise, we first need to
define a appropriate compatibility condition between a trivialized pencil bundle and a metric.

\begin{definition} ---
A metric on a manifold $M$ is called \emph{compatible} with an indexed pencil bundle if the angle between two
given plane bundles $\xi_{\theta_1}$ and $\xi_{\theta_2}$ is constant on $M$.

A moving frame $(X_1,\,X_2,\,X_3)$ is said to be \emph{adapted} to an indexed pencil bundle if $X_3$ is parallel
to the axis of the pencil and if there are forms $\omega_1$ and $\omega_2$ in the defining family of 1-forms with
$\omega_1(X_2)=\omega_2(X_1)=0$ everywhere on $M$.
\end{definition}

\noindent {\bf Remark} --- A metric for which a moving frame $(X_1,\,X_2,\,X_3)$ adapted to some indexed pencil
bundle is orthonormal is compatible with the indexed pencil bundle whose trivialization is induced by the family
$\mathcal{S}^1\{\omega_1,\,\omega_2\}$, with $\omega_1(X_2)=\omega_2(X_1)=0$. \vspace{0.4cm}

The rotation speed of a vector field about some other vector field along a given curve is defined as follows~:

\begin{definition} ---
Let $(M,g)$ be an oriented Riemannian manifold, $Y$ and $Z$ unit vector fields and $\gamma$ a curve on $M$. We
will call the quantity
$$ R_\gamma(Y,Z)=g(\nabla_{\dot{\gamma}}Y, Z \times Y)$$ the \emph{rotation speed} of Y about Z along $\gamma$
with respect to the parallel transport.
\end{definition}

We choose this terminology, because $g(\nabla_{\dot{\gamma}}Y, Z \times Y)$ can be seen as the component of the
projection of the covariant derivative of $Y$ along $\gamma$ onto the plane orthogonal to $Z$ which is orthogonal
to the projection of $Y$ onto this same plane.

\begin{definition} ---
Let $\mathcal{F}$ be an indexed pencil bundle on an oriented Riemannian manifold $(M,g)$. The \emph{spinning
direction} of $\mathcal{F}$ is defined to be an orientation of the axis bundle such that the rotation speed of
$X_1$ about $X_3$ along every integral curve of this axis bundle is non-negative for a positively oriented adapted
moving frame $(X_1,X_2,X_3)$ such that $X_3$ is positively oriented on the axis bundle.
\end{definition}

The following definition has been modeled on the properties of the fundamental example on $\sph^3$ developed
below (see Examples \ref{ex sphere} and \ref{ex contact tops}).

\begin{definition} \label{toupie} ---
Let $(M,g)$ be an oriented Riemannian manifold of dimension~3. An indexed pencil bundle on $M$ is called a
\emph{top} if there is an orthonormal moving frame $(X_1,\,X_2,\,X_3)$ adapted to the pencil bundle such that
\begin{enumerate}
\item[(i)] Along any geodesic $\gamma$, the angle between $\dot{\gamma}$
and $X_3$ is constant. \vspace{-2mm}

\item[(ii)] Along any geodesic $\gamma$ which is transverse to the axis bundle, the rotation speed of $X_3$ about
$\dot{\gamma}$ is constant with respect to the parallel transport. \vspace{-2mm}

\item[(iii)] Along any geodesic $\gamma$, the rotation speed of $X_1$ about
$X_3$ is constant with respect to the parallel transport. \vspace{-2mm}

\item[(iv)] Along any pair of geodesics $\gamma_1$ and $\gamma_2$ such that the angles between $\dot{\gamma_1}$ and
$X_3$ and between $\dot{\gamma_2}$ and $X_3$ are equal, the rotation speed of $X_1$ about $X_3$ is the same. \vspace{-2mm}

\end{enumerate}
All geodesics are supposed to be parameterized by arc length.

If every plane bundle defined by the trivialization is a contact structure, the top will be called a
\emph{contact top}; if all these plane bundles are integrable, it will be called an \emph{integrable top}.
\end{definition}

\noindent {\bf Remarks} ---
\begin{enumerate}
\item A top is completely determined by a moving frame satisfying the four conditions above.
\item If a top is determined by an orthonormal moving frame $(X_1,X_2,X_3)$, the underlying indexed pencil bundle is
defined by the dual forms $\omega_1$ and $\omega_2$ which correspond to $X_1$ and $X_2$.
\item Condition (i) implies that the integral curves of the axis bundle are geodesics (see also Lemma \ref{lem D}).
\end{enumerate}

We will see that a top is always either a contact top or an integrable top (see Corollary \ref{type of tops}).
Contact tops are always defined by taut contact circles (see the definitions in section \ref{def contact}. Some
examples will also be given there and in section \ref{integrable tops}.)

The first main result of this paper is a characterization of tops through the properties of the Lie brackets of
an adapted orthonormal moving frame.

\begin{theorem} --- \label{car toupies}

An orthonormal global moving frame $(X_1,\,X_2,\,X_3)$ on an oriented Riemannian manifold $M$ determines a top if
and only if the corresponding Lie brackets are of the following type~:
\begin{eqnarray} \label{crochets}
\left \{ \begin{array}{lcl}
\left[ X_1,X_2\right]&=& c \, X_3 \\
\left[ X_2,X_3\right]&=& k \, X_1 \\
\left[ X_3,X_1\right]&=& k \, X_2,
\end{array} \right .
\end{eqnarray}
where $c$ and $k$ are real numbers.
\end{theorem}

As a corollary, we obtain a description of the manifolds which admit tops.

\begin{cor} \label{var toupies} ---
An oriented and complete connected manifold $M$ carries a top $\mathcal{T}$ if and only if $M$ is diffeomorphic
to the left-quotient of one of the following Lie groups by a discrete subgroup~:

$\left. \begin{array}{l}
\widetilde{E}_2 \\
SU(2) \\
\widetilde{SL}(2,\R)
\end{array} \right\}$ for a contact top;
$\quad \left. \begin{array}{l}
\R^3 \\
Nil^3
\end{array} \right\}$ for an integrable top,

\noindent where $\widetilde{E}_2$ is the universal cover of the group of direct isometries of the plane and
$Nil^3$ the $3$-dimensional Heisenberg group.
\end{cor}

We are now interested in describing the metrics on these manifolds for which tops can be constructed. They have
to be adapted to the geometric properties of tops.

\begin{definition} --- \label{spinning metric}
A metric $g$ on a 3-manifold $M$ is called a \emph{spinning metric} if
\begin{itemize}
\item there is a unit vector field $Z$ on $M$,
called \emph{pivot field}, which is geodesic and Killing, such that the sectional curvature of any plane only
depends on its angle with $Z$ and if
\item the extremal sectional curvatures are constant on $M$.
\end{itemize}
\end{definition}

These are the only metrics to which tops can be associated, as the following theorem states.

\begin{theorem} \label{spin thm} ---
Let $(M,g)$ be an oriented Riemannian $3$-manifold. If $(M,g)$ admits a top, then $g$ is a spinning metric. If
$g$ is a spinning metric and $H_1(M,\R)=0$, then there exists a spinning metric on $(M,g)$.
\end{theorem}

Sections \ref{def contact} and \ref{integrable tops} contain background information and examples concerning
contact tops and integrable tops. In section \ref{bundles and tops}, we prove Theorem \ref{car toupies} and
several corollaries. In section \ref{spinning metrics}, we analyze the properties of spinning metrics and prove
\mbox{Theorem \ref{spin thm}}. In section \ref{compatible metrics}, we determine the metrics for which a given
contact pencil bundle defines a top. In section \ref{classification}, some partial classification results are
discussed, especially the uniqueness problem of a top for a given spinning metric.

In the last part of this paper, a connection to Sasakian geometry is made. In particular, we discuss which
Sasakian structures define tops and which tops give rise to a Sasakian structure. \vspace{4mm}

To conclude, let us remark that it is possible to consider tops in higher dimensions, but there are considerable
additional difficulties. The first one is the problem to find a good definition, because the common subspace
bundle of a family of hyperplanes in $\mathbb{R}^n$ can be of different dimensions, depending on the number of
generators of the defining family of 1-forms, and as soon as the codimension of the common subspace bundle is
greater than 2, one can not use the concept of rotation speed any more. Moreover, the study of contact
$p$-spheres on higher-dimensional manifolds is much less developed than for contact circles on 3-manifolds.

\section{Contact circles and contact tops} \label{def contact}

Some examples coming from contact geometry have been an important motivation to introduce contact tops. The study
of contact circles has been initiated by H. Geiges and J. Gonzalo in 1995 (see \cite{GG1}).

\begin{definition} ---
If all non-trivial, normalized linear combinations of two contact forms are contact forms, this family is called a
\emph{contact circle}. We note $\mathcal{S}^1_c\{\omega_1,\,\omega_2\}:=\{ \lambda_1 \, \omega_1 + \lambda_2 \,
\omega_2, \lambda_1^2 + \lambda_2^2 = 1 \}$.

Similarly, a contact sphere is generated by three contact forms.
\end{definition}

H. Geiges and J. Gonzalo distinguish a certain class of contact circles, defined as follows~:
\begin{definition} ---
A contact circle $\mathcal{S}^1_c\{\omega_1,\,\omega_2\}$ on a 3-manifold $M$ is said to be \emph{taut} if all
its elements define the same volume form, that is, if $\omega \wedge d\omega$ is constant on
$\mathcal{S}^1_c\{\omega_1,\,\omega_2\}$.
\end{definition}

\begin{example} --- \emph{Contact circle on $\tor^3$} \label{ex tore}

On the 3-torus with pseudo-coordinates $(\theta_1,\,\theta_2,\,\theta_3)$, the forms
\begin{eqnarray*}
\omega_1&=&\phantom{-}\cos(n\theta_1) \, d\theta_2+\sin(n\theta_1) \, d\theta_3\\
\omega_2&=&-\sin(n\theta_1) \, d\theta_2+\cos(n\theta_1) \, d\theta_3
\end{eqnarray*}
generate a contact circle, for $n \in \mathbb{N}^*$.

Indeed, for $\omega=\lambda_1\,\omega_1+\lambda_2\,\omega_2$, with $\lambda_1^2+\lambda_2^2=1$, we have $ \omega
\wedge d\omega=-n\,d\theta_1 \wedge d\theta_2 \wedge d\theta_3, $ so all non-trivial linear combinations are
contact forms. They all define the same volume form, so the contact circle is taut.
\end{example}

\begin{example} --- \emph{Contact sphere on $\sph^3$} \label{ex sphere}

We consider the 3-sphere as the unit sphere of the quaternionic space $\ha$. On this sphere, we have 3 independent
contact forms, induced by the following forms on $\ha$~:
\begin{eqnarray*}
\alpha_q&=&<qi,\,dq>=q_1\,dq_2-q_2\,dq_1+q_4\,dq_3-q_3\,dq_4,\\[1mm]
\beta_q&=&<qj,\,dq>=q_1\, dq_3-q_3\,dq_1+q_2\,dq_4-q_4\,dq_2,\\[1mm]
\gamma_q&=&<qk,\,dq>=q_3\, dq_2-q_2\,dq_3+q_1\,dq_4-q_4\,dq_1.
\end{eqnarray*}

The induced forms generate a contact sphere. Indeed, any normalized linear combination $ \omega:=\lambda_1 \,
\alpha + \lambda_2 \, \beta + \lambda_3 \, \gamma $ satisfies~:
$$ \omega \wedge d\omega \wedge\,(q_1\,dq_1+q_2\,dq_2+q_3\,dq_3+q_4\,dq_4)=dq_1\wedge \,dq_2\wedge \,dq_3\wedge \,dq_4,$$
which is a volume form on $\ha$. This volume form is independent of the parameters
$\lambda_1,\lambda_2,\lambda_3$, so $\alpha$, $\beta$ and $\gamma$ induce a taut contact sphere on $\sph^3$.
\end{example}

We have another regularity condition for contact circles, with a more obvious geometric meaning than tautness~:
\begin{definition} ---
A contact circle $\mathcal{S}^1_c\{\omega_1,\,\omega_2\}$ is said to be \emph{round} if the Reeb vector field of
any element $\omega=\lambda_1 \, \omega_1 + \lambda_2 \, \omega_2$ of the contact circle is given by $R=\lambda_1
R_1 + \lambda_2 R_2$, where $R_1$ and $R_2$ are the Reeb vector fields of $\omega_1$ and $\omega_2$.
\end{definition}

On 3-manifolds, roundness and tautness are equivalent, but on higher-dimensional manifolds, these properties are
independent (see \cite{Z}). \vspace{4mm}

The first important question is about the existence of contact circles. For closed manifolds, H. Geiges and J.
Gonzalo give a general answer~: \vspace{4mm}

\noindent {\bf Theorem} (see \cite{GG2}) --- \emph{Every closed and orientable 3-manifold admits a contact
circle.} \vspace{1mm}

For taut contact circles, which have richer geometric properties, the existence problem is more complex~:
\vspace{4mm}

\noindent {\bf Theorem} (see \cite{GG1}) --- \textit{A closed and orientable 3-manifold $M$ admits a taut contact
circle if and only if $M$ is diffeomorphic to the left-quotient of a Lie group $G$ by the action of a discrete
subgroup, where $G$ is either $SU(2)$ or $\widetilde{SL}(2,\R)$, the universal cover of $SL(2,\R)$, or
$\widetilde{E}_2$, the universal cover of the Euclidean group.} \vspace{4mm}

A contact circle defines a contact pencil bundle in an obvious way. If we write the contact circle as
$\{\omega^\theta=\cos\theta \,\omega_1+\sin\theta \,\omega_2, \theta\in\left[0,2\pi\right[\}$, the corresponding
trivialization is given by $\tau(\ker(\omega^\theta))=\theta \mod \pi$.

A round contact circle defines a privileged plane bundle, which is transverse to the axis bundle
$\ker\omega_1\cap\ker\omega_2$. It is spanned by the Reeb vector fields $R_1$ and $R_2$ of any two generating
forms $\omega_1$ and $\omega_2$.

\begin{example} \label{ex contact tops} --- The contact circles in Example \ref{ex tore} and the contact circles
generated by any two of the three forms of Example \ref{ex sphere} define contact tops, for the flat metric on
$\tor^3$ and the usual metric
on $\sph^3$.

To see this in the case of the 3-sphere, consider the standard metric and an orthonormal frame $(X_1,X_2,X_3)$
dual to the induced contact forms $(\tilde{\alpha}, \tilde{\beta}, \tilde{\gamma})$. This frame satisfies
$$ \left[ X_1,X_2 \right] = 2 X_3, \quad \left[ X_2,X_3 \right] = 2 X_1 \quad \textrm{and} \quad \left[ X_3,X_1
\right] = 2 X_2.$$

On $\sph^3$, a vector field $Z=a_1 X_1 + a_2 X_2 + a_3 X_3$ is geodesic if and only if its coefficients are
constant on $\sph^3$. So $X_3$ is a geodesic vector field. Let $\rho$ with $\dot{\rho}=a_1 X_1 + a_2 X_2 + a_3
X_3$ be a geodesic.

Then the rotation speed of $X_1$ about $X_3$ along $\rho$ is
$$ R_\rho(X_1,X_3) = g(\nabla_{\dot{\rho}}X_1,X_2) = 2 a_3, $$
so it is constant and only depends on the angle between $X_3$ and $\dot{\rho}$, that is, on $a_3$, which is
constant.

The rotation speed of $X_3$ about $\dot{\rho}$ along $\rho$ is
$$ R_\rho(X_3,\dot{\rho}) = a_1^2 + a_2^2, $$
so it is constant, too.

Thus, the contact circle $\mathcal{S}_c^1\{\tilde{\alpha},\tilde{\beta}\}$ defines a contact top.
\end{example}

\section{Integrable tops} \label{integrable tops}
An integrable top on a Riemannian manifold $M$ defines a circle of foliations of $M$ which are everywhere
transverse to each other. A trivial example is the following~:

\begin{example} \label{ex int top} ---
A pencil of planes at the origin of $\R^3$ which is parallel transported to all other points of $\R^3$ defines an
integrable top with respect to the Euclidean metric. In this example, all rotations speeds are zero and all
angles are constant, so the conditions of Definition \ref{toupie} are trivially satisfied. All plane bundles
defined by this top are integrable.
\end{example}

Another example is the following indexed pencil bundle on the Heisenberg group.

\begin{example} --- \label{Heisenberg group}
We represent the Heisenberg group as $\R^3$ with the group operation $(x,y,z)*(x',y',z')=(x+x',y+y',z+z'+x\,y')$,
and with the metric for which the system $(X_1,X_2,X_3)=(\frac{\partial}{\partial x}\,, \frac{\partial}{\partial
z}+x\,\frac{\partial}{\partial y}\,, \frac{\partial}{\partial y})$ is orthonormal. The tangent vector field
$X=a_1\,X_1+a_2\,X_2+a_3\,X_3$ of a geodesic $\gamma$ in this space satisfies $\nabla_XX=0$, thus
$$X(a_1)=0,\quad X(a_2)=0,\quad X(a_3)=0.$$

So a geodesic $\gamma$ with unit speed is defined by $ \dot{\gamma}(t)=a_1X_1 + a_2X_2 + a_3 X_3, $ for
$a_1,a_2,a_3 \in \R$.

Thus, the vector field $X_3$ is geodesic and the rotation speed of $X_1$ about $X_3$ along $\gamma$ is
$$ R_\gamma(X_1,X_3)=-\frac{a_3}{2}\,, $$
so it is constant and it only depends on the angle between $X_3$ and $\dot{\gamma}$, that is, on $a_3$, which is
constant, too.

The rotation speed of $X_3$ about $\dot{\gamma}$ along $\gamma$ is
$$ R_\gamma(X_3,\dot{\gamma})=\frac{1}{2}\,(a_1^2+a_2^2), $$
so it is constant, as $a_1^2+a_2^2=1-a_3^2$.

Thus, $(X_1, X_2, X_3)$ defines a top on the Heisenberg group.
\end{example}

\section{Indexed pencil bundles and tops} \label{bundles and tops}
We are now going to have a closer look at indexed plane pencil bundles.

Definition \ref{toupie} describes tops in a geometric way. We will now characterize them using Lie brackets, in
order to have some additional working tools. This is done by Theorem \ref{car toupies}, which implicitly gives us
much important information about the manifolds which carry tops and about the plane bundles which are associated
to tops, as a number of corollaries will show.

We first need a lemma which relates several geometric properties of a vector field on a Riemannian manifold.

\begin{lemma} \label{lem D} ---
Let $Z$ be a unitary vector field on a Riemannian manifold. Then the following conditions are equivalent:
\begin{enumerate}
\item[(i)] $Z$ is a Killing vector field and its integral curves are geodesics.
\item[(ii)] Along any geodesic $\gamma$, the angle between $\dot{\gamma}$ and $Z$ is constant.
\item[(iii)] For every orthonormal frame $(X_1,X_2,Z)$, there is a real-valued function $k$, such that
$\left[Z,X_1\right]=kX_2$ and $\left[X_2,Z\right]=kX_1$.
\end{enumerate}
\end{lemma}

\noindent {\bf Proof~:} {\it (i)$\Leftrightarrow$(iii)} : Let $(X_1,X_2,Z=X_3)$ be an orthonormal frame on $M$.
We write $$\left[ X_i,X_j\right]=\sum_{i\neq j} c_{ij}^k X_k.$$ The integral curves of $Z$ are geodesics if and
only if $\nabla_Z Z=0$, which is equivalent to $c_{31}^3=c_{23}^3=0$. For all calculations of this type, we write
$\Gamma_{ij}^k=g(\nabla_{X_i}X_j, X_k)$ and we use the relation
$\Gamma_{ij}^k=\frac{1}{2}\,(c_{ij}^k+c_{ki}^j+c_{kj}^i)$ (see e.g. \cite{Hel}, p. 48).

$Z$ is a Killing vector field if and only if for every pair of vector fields ($X$,$Y$) on $M$, we have
$g(\nabla_X Z,Y)+g(X,\nabla_Y Z)=0$, which means that $c_{31}^1=c_{23}^2=0$ (for $X=Y=X_1$ and for $X=Y=X_2$) and
$c_{23}^1=c_{31}^2$ (for $X=X_1$ and $Y=X_2$).

\vspace{4mm} {\it (ii)$\Leftrightarrow$(iii)} : For a geodesic $\gamma$ with
$\dot{\gamma}(s)=a_1(s)\,X_1+a_2(s)\,X_2+a_3(s)\,X_3$, we have
$$ \frac{\partial}{\partial s}(g(X_3,\dot{\gamma}))=g(\nabla_{\dot{\gamma}}X_3, \dot{\gamma})=a_1^2 \,c_{13}^1 +
a_2^2\,c_{23}^2 + a_1\,a_2(c_{13}^2 + c_{23}^1) + a_1\,a_3\,c_{13}^3 + a_2\,a_3\,c_{23}^3. $$
So the constancy of the angle between $Z$ and $\dot{\gamma}$ along any geodesic $\gamma$ is equivalent to the
property {\it (iii)}. \finthm \vspace{4mm}

\noindent {\bf Remark} --- It follows from this Lemma that a unitary vector field which defines the axis bundle
of a top is a Killing vector field. \vspace{4mm}

\noindent {\bf Proof} of Theorem \ref{car toupies}~: Let $(X_1,\,X_2,\,X_3)$ be a positively oriented orthonormal
moving frame which defines a top. Again, we write the corresponding Lie brackets as
$$\left[ X_i,X_j\right]=\sum_{i\neq j} c_{ij}^k X_k.$$

Furthermore, let $\gamma$ be a geodesic such that $\dot{\gamma}$ has unit length. We write \mbox{$\dot{\gamma}=a_1
X_1 + a_2 X_2 + a_3 X_3$}, where the $a_i$ are functions on $M$. We will express the geometric conditions of
Definition \ref{toupie} in terms of the functions $c_{ij}^k$~:
\begin{enumerate}
\item[{\it (i)}] By Lemma \ref{lem D}, the property that the angle between $X_3$ and $\dot{\gamma}$ is constant
along any geodesic $\gamma$ is equivalent to the identities $c_{31}^1=c_{23}^2=c_{31}^3=c_{23}^3=0$ and $c_{23}^1=c_{31}^2$
on $M$. In that case, it also follows that $a_3$ is constant along $\gamma$.

\item[{\it (iv)}] Under the assumption that property {\it (i)} is satisfied, we have
\begin{eqnarray*}
R_\gamma(X_1,X_3) &=& a_1 \Gamma_{11}^2 + a_2 \Gamma_{21}^2 + a_3 \Gamma_{31}^2 \\
&=& a_1 c_{21}^1 + a_2 c_{21}^2 + a_3 (c_{23}^1-\frac{1}{2}c_{12}^3).
\end{eqnarray*}
Thus, the rotation speed of $X_1$ about $X_3$ along $\gamma$ does not depend on the angles between $\dot{\gamma}$
and $X_1$ or $X_2$, if and only if $ c_{21}^1=c_{21}^2=0.$

\item[{\it (ii)}] Suppose now that $\gamma$ is transverse to the axis bundle. If properties {\it (i)} and {\it (iv)}
are satisfied, we have
\begin{eqnarray*}
R_\gamma(X_3,\dot{\gamma}) &=& a_2\, g(\nabla_{\dot{\gamma}} X_3,X_1) - a_1\, g(\nabla_{\dot{\gamma}} X_3,X_2) \\
&=& c_{12}^3 (a_1^2 + a_2^2).
\end{eqnarray*}
Thus, the constancy of the rotation speed of $X_3$ about $\dot{\gamma}$ along $\gamma$ is equivalent to the
constancy of $c_{12}^3$ along $\gamma$. Indeed, $a_1^2 + a_2^2$ is constant along $\gamma$, since $a_3$ is, by
{\it (i)}. As this is true for any choice of $\gamma$ (transverse to $X_3$), it follows that $c_{12}^3$ is
constant on $M$ if and only if the rotation speed of $X_3$ about $\dot{\gamma}$ is constant along any geodesic
$\gamma$.

\item[{\it (iii)}] Now we assume that properties {\it (i)}, {\it (ii)} and {\it (iv)} are satisfied. So we have
$$ R_\gamma(X_1,X_3) = a_3 (c_{23}^1-\frac{1}{2}\,c_{12}^3). $$
Thus, the rotation speed of $X_1$ about $X_3$ is constant along $\gamma$ if and only if $c_{23}^1$ is constant
along $\gamma$. Again, it follows that $c_{23}^1$ is constant on $M$ if and only if the rotation speed of $X_1$
about $X_3$ is constant along any geodesic $\gamma$. \finthm
\end{enumerate}

Theorem \ref{car toupies} allows us to describe the nature of the indexed plane bundles associated to a top as
well as the manifolds where tops can be constructed.

We first prove Corollary \ref{var toupies}, stated in the introduction. \vspace{4mm}

\noindent {\bf Proof} of  Corollary \ref{var toupies}~: Let $(X_1,\,X_2,\,X_3)$ be a moving frame on $M$ which
defines $\mathcal{T}$ and so satisfies the relations (\ref{crochets}). This defines on $M$ a locally free action
of the simply connected Lie group $\mathcal{G}$ whose Lie algebra is defined by these relations. According to the
classification of Lie algebras of dimension 3, these Lie algebras are $\mathfrak{so}(3,\R)$ (for $ck>0$),
$\mathfrak{sl}(2,\R)$ (for $ck<0$), the Lie algebra associated to $E_2$ (for $c=0,\,k \neq 0$), the Lie algebra
associated to $Nil^3$ (for $c \neq 0, \, k=0$), or $\R^3$ (for $c=k=0$). Any orbit of this action is
diffeomorphic to a quotient of $\mathcal{G}$ by the stabilizer, which is a discrete subgroup. So every orbit is
an open submanifold of $M$ and its complement is a union of orbits, so it is open also. Hence, the orbit is a
connected component of $M$ and the action is transitive, $M$ being connected. Thus, $M$ is diffeomorphic to a
quotient of one of the corresponding Lie groups by a discrete subgroup given by the isotropy subgroup of a point
of $M$.

On the other hand, given one of these Lie groups, which we call $\mathcal{G}$, and a discrete subgroup
$\mathcal{S}$, let $(X_1,\,X_2,\,X_3)$ be a global moving frame of $\mathcal{G}$, where each $X_i$ is a
left-invariant vector field coming from a generator of the associated Lie algebra, such that (\ref{crochets}) is
satisfied for some constants $c$ and $k$. In particular, this moving frame is invariant under the left-action of
$\mathcal{S}$. Hence, it induces a moving frame on the quotient $\mathcal{S} \backslash \mathcal{G}$ which
satisfies the relations (\ref{crochets}). According to Theorem \ref{car toupies}, there exists a top on
$\mathcal{S} \backslash \mathcal{G}$. \finthm

\begin{cor} --- \label{type of tops}
All plane bundles defined by a given top are of the same nature. They are either integrable or contact structures.
\end{cor}

\noindent {\bf Proof :} Let $\omega_1,\,\omega_2$ and $\omega_3$ be the dual forms associated to the vector fields
$X_1,\,X_2$ and $X_3$ of an orthonormal moving frame of $M$ satisfying the relations (\ref{crochets}), that is,
\begin{eqnarray} \label{formes}
\left \{ \begin{array}{lcl}
d\omega_3 &=& -c \, \omega_1 \wedge \omega_2 \\
d\omega_2 &=& -k \, \omega_3 \wedge \omega_1 \\
d\omega_1 &=& -k \, \omega_2 \wedge \omega_3. \\
\end{array} \right .
\end{eqnarray}

Then $\omega_1$ and $\omega_2$ are the generating forms of the top which is determined by $(X_1,X_2,X_3)$. The
plane bundles defined by the trivialization of this top are the kernels of the forms $\omega^\theta=\cos \theta \,
\omega_1 + \sin \theta \, \omega_2$.

If $k=0$, then $\omega_1$ and $\omega_2$ are integrable, as well as their linear combinations, so the pencil
bundle is integrable.

If $k\neq 0$, then $\omega_1$ and $\omega_2$ are contact forms, as well as their linear combinations. So these
two forms generate a contact circle $\mathcal{S}^1_c\{\omega_1,\,\omega_2\}$ and the corresponding pencil bundle
is a contact pencil bundle. \finthm

\begin{definition} (see \cite{GG3}) ---
A contact circle $\mathcal{S}^1_c\{\omega_1,\,\omega_2\}$ is called a K-Cartan structure if there is a 1-form
$\omega_3$ such that
$$ d\omega_1=\omega_2 \wedge \omega_3, \quad d\omega_2=\omega_3 \wedge \omega_1 \quad \textrm{and} \quad d\omega_3=K
\omega_1 \wedge \omega_2,$$ for some real number $k$. In particular, a K-Cartan structure is taut.
\end{definition}

\begin{cor} \label{cor str Cartan} ---
Any contact top is defined by a K-Cartan structure, for \mbox{$ K \in \{ -1,0,1 \} $}.
\end{cor}

\noindent {\bf Proof~:} In the proof of the preceding Corollary, we show that a contact top is always defined by
a taut contact circle. This contact circle is in fact a K-Cartan structure, by (\ref{formes}), where K can be
chosen as 0, 1 or -1. \finthm

These results about contact tops can be related to Geiges' and Gonzalo's theorem about the classification of
closed and orientable 3-manifolds which carry taut contact circles (see \cite{GG1} or section \ref{introduction}),
in the following way.

Corollary \ref{var toupies} implies that on all manifolds listed in this theorem taut contact circles exist,
since any contact top is defined by a taut contact circle (see Corollary \ref{cor str Cartan}).

On the other hand, Proposition 3.5 of \cite{GG3} states that every conformal class of taut contact circles on a
compact left-quotient of $\widetilde{SL}_2$ (resp. $\widetilde{E}_2$) contains a K-Cartan structure, for $K=-1$
resp. 0. Such a conformal class is defined as the set of multiples of a given contact circle, that is, where both
generators are multiplied by the same positive function. This does not change the associated contact structures,
so a conformal class of taut contact circles defines a trivialized plane pencil bundle, which is a top if this
conformal class contains a K-Cartan structure. Thus, on these manifolds, every taut contact circle defines a
contact top. \vspace{4mm}

In the following sections, we are going to have a closer look on the metrics which can be associated to tops. The
problems are the following : {\it For which metrics can tops be constructed? How far does a given metric or a
given top determine the tops and metrics we can associate to it?} We will partially answer these questions in the
following sections.

\section{Spinning metrics} \label{spinning metrics}
In the relations (\ref{crochets}), $X_1$ and $X_2$ play symmetric parts, so a metric for which a moving frame
satisfying (\ref{crochets}) is orthonormal has to be, in some sense, homogeneous around the vector field $X_3$,
that is, invariant by rotation about the axis given by $X_3$.

Let us see some examples of spinning metrics (see Definition \ref{spinning metric})~:
\begin{enumerate}
\item[{\it i)}] A metric with curvature zero is a spinning metric, if there is a geodesic Killing vector field.
An example is given by the flat metric on $\tor^3$. \vspace{-2mm}

\item[{\it ii)}] Let $(Y_1,\,Y_2,\,Y_3)$ be a moving frame on $\sph^3$ which defines the usual metric and where
$Y_3$ is a geodesic Killing vector field. Then any metric for which the moving frame $(a\,Y_1,\,a\,Y_2,\,b\,Y_3)$,
with $a,\,b \in \R^*$, is orthonormal, is a spinning metric. \vspace{-2mm}

\item[{\it iii)}] Another interesting case is the situation where the sectional curvature of the planes
which contain $Z$ is positive and where the plane orthogonal to $Z$ has negative sectional curvature. An example
of this type is given on the Heisenberg group (see Example \ref{Heisenberg group}).
\end{enumerate}

\noindent {\bf Remark} --- Metrics with negative constant curvature are not spinning metrics, because they do not
admit non-trivial Killing vector fields. According to the calculations in the proof of the following lemma, any
plane which contains a pivot field has non-negative sectional curvature.

\begin{lemma}\label{crochets metr piv} ---
Let $g$ be a spinning metric on a manifold $M$ and let $Z$ be a pivot field of $g$. If $(X_1,\,X_2,\,Z)$ is an
orthonormal moving frame, the Lie brackets of these vector fields are of the following type~:
\begin{eqnarray}\label{crochets piv}
\left \{ \begin{array}{lcl}
\left[ X_1,X_2\right]&=&c_{12}^1 X_1+c_{12}^2 X_2+c_{12}^3 Z \\
\left[ X_2,Z\right]&=&c_{23}^1 X_1 \\
\left[ Z,X_1\right]&=&c_{23}^1 X_2,
\end{array} \right .
\end{eqnarray}
where $c_{12}^3$ is a constant.
\end{lemma}

\noindent {\bf Proof~:} Lemma \ref{lem D} implies that the Lie brackets satisfy (\ref{crochets piv}) and it only
remains to show that $c_{12}^3$ is a constant.

In the following, we will use the curvature identities (see for example \cite{Hel}, p.69)
\begin{eqnarray*}
g(R(X,Y)Z,T) &=& g(R(Z,T)X,Y) \\
g(R(X,Y)Z,T) &=& -g(R(Y,X)Z,T) \\
g(R(X,Y)Z,T) &=& -g(R(X,Y)T,Z),
\end{eqnarray*}
where $R(X,Y)=\nabla_X \nabla_Y - \nabla_Y \nabla_X - \nabla_{\left[X,Y\right]}$, and the relation (see
\cite{Hel}, p.45)
$$ g(R(X_i,X_j)X_k,X_s)=\sum_l (\Gamma_{jk}^l \Gamma_{il}^s - \Gamma_{ik}^l \Gamma_{jl}^s) +
X_i(\Gamma_{jk}^s) - X_j(\Gamma_{ik}^s) - \sum_l c_{ij}^l \Gamma_{lk}^s.$$ A plane spanned by two vectors $$
X=\cos \theta X_1 + \sin \theta X_2 \quad \textrm{and} \quad Y=- \sin \varphi \sin \theta X_1 + \sin \varphi \cos
\theta X_2 + \cos \varphi Z $$ has sectional curvature
\begin{eqnarray} \label{spin metric}
    && - g(R(X,Y)X,Y) \nonumber \\
    &=& - \sin^2 \varphi \, g(R(X_1,X_2)X_1,X_2) \nonumber \\
    && - \cos^2 \varphi \, \Big(\cos^2 \theta \, g(R(X_1,Z)X_1,Z) + \sin^2 \theta \, g(R(X_2,Z)X_2,Z) \Big)
    \nonumber \\
    && - 2 \cos \theta \, \sin \varphi \, \cos \varphi \, g(R(X_1,X_2)X_1,Z) \nonumber \\
    && - 2 \sin \theta \, \sin \varphi \, \cos \varphi \, g(R(X_1,X_2)X_2,Z) \nonumber \\
    && - 2 \cos \theta \, \sin \theta \, \cos^2 \varphi \, g(R(X_1,Z)X_2,Z) \nonumber \\
    &=& \sin^2 \varphi \, \Big( - \frac{3}{4}(c_{12}^3)^2 - (c_{12}^1)^2 - (c_{12}^2)^2 + c_{12}^3 \, c_{23}^1 -
    X_2(c_{12}^1) + X_1(c_{12}^2) \Big) \\
    && + \frac{1}{4} \, \cos^2 \varphi \, (c_{12}^3)^2 \nonumber \\
    && + \cos \theta \, \sin \varphi \, \cos \varphi \, X_1(c_{12}^3) \nonumber \\
    && + \sin \theta \, \sin \varphi \, \cos \varphi \, X_2(c_{12}^3) \nonumber \\
    && + \cos \theta \, \sin \theta \, \cos^2 \varphi \, Z(c_{12}^3) \nonumber.
\end{eqnarray}
As $g$ is a spinning metric, this expression does not depend on $\theta$, hence
\begin{equation*}
    X_1(c_{12}^3) = 0, \quad\;
    X_2(c_{12}^3) = 0, \quad\;
    Z(c_{12}^3) = 0.
\end{equation*}
Thus, $c_{12}^3$ is constant on $M$. \finthm

\begin{lemma} \label{courbures pivotantes} ---
Let $g$ be a spinning metric on a $3$-manifold $M$. Let $Z$ be a pivot field of $g$, let $\alpha$ be the sectional
curvature of any plane containing $Z$ and let $\beta$ be the sectional curvature of the plane which is orthogonal
to $Z$. Then $\alpha$ and $\beta$ are the extremal sectional curvatures of $g$ and the sectional curvature of a
plane which forms an angle $\varphi$ with $Z$ is
  $$\alpha \, \cos^2 \varphi+\beta \, \sin^2 \varphi.$$
\end{lemma}

\noindent {\bf Proof~:} This is another consequence of the calculation of the sectional curvature of an arbitrary
plane for a spinning metric in the proof of Lemma \ref{crochets metr piv}. \finthm

From now on, when we talk about the extremal sectional curvatures $\alpha$ and $\beta$ of a spinning metric,
$\alpha$ will be the sectional curvature of any plane which contains a pivot field $Z$, and $\beta$ will be the
sectional curvature of the plane orthogonal to $Z$. This does not depend on the choice of $Z$, by Lemma \ref{Z
unique} below. \vspace{4mm}

Theorem \ref{spin thm}, stated in the introduction, specifies the close relationship between tops and spinning
metrics.

To prove it, we need the following lemma~:
\begin{lemma} \label{Z unique} ---
If $g$ is a spinning metric on a $3$-manifold $M$ and if $(M,g)$ is not a space of constant curvature, then the
pivot field of $g$ is unique up to the sign.
\end{lemma}

\noindent {\bf Proof~:} This is a consequence of Lemma \ref{courbures pivotantes}. \finthm

\noindent {\bf Proof} of Theorem \ref{spin thm}~: Let us see why $g$ is a spinning metric if there is a top
$\mathcal{T}$ on $(M,g)$. Let $\mathcal{T}$ be defined by an orthonormal moving frame $(X_1, X_2, X_3)$ satisfying
(\ref{crochets}). So $X_3$ is geodesic and Killing. According to the calculations in the proof of Lemma
\ref{crochets metr piv}, the sectional curvature of the plane spanned by the vectors
$$ X=\cos \theta X_1 + \sin \theta X_2, \quad \textrm{and} \quad Y=- \sin \varphi \sin \theta X_1 + \sin \varphi \cos
\theta X_2 + \cos \varphi X_3,$$
for arbitrary angles $\theta$ and $\varphi$, is $\cos^2 \varphi \, \frac{c^2}{4}
+ \sin^2\varphi \, (ck - \frac{3}{4} \, c^2)$. This value does not depend on $\theta$. The extremal sectional
curvatures are
\begin{equation} \label{courb extr}
\alpha = \frac{c^2}{4} \quad \textrm{and} \quad \beta = ck - \frac{3}{4} \, c^2,
\end{equation}
so they are constant. Hence, $g$ is a spinning metric. \vspace{0.4cm}

Let now $g$ be a spinning metric on $M$ and $Z$ an associated pivot field and let $(X_1,\,X_2,\,Z)$ be a
positively oriented orthonormal moving frame. By Lemma \ref{crochets metr piv}, these vector fields satisfy
(\ref{crochets piv}), where $c_{12}^3$ is a constant.

The Jacobi condition yields
\begin{eqnarray}\label{Jacobi2}
\left\{ \begin{array}{rcl}
    X_1(c_{23}^1) - c_{12}^2 \, c_{23}^1 + Z(c_{12}^1) &=& 0 \\
    X_2(c_{23}^1) + c_{12}^1 \, c_{23}^1 + Z(c_{12}^2) &=& 0 \\
    Z(c_{12}^3) &=& 0.
\end{array} \right.
\end{eqnarray}
Moreover, the extremal sectional curvatures of a spinning metric are constant, so
\begin{equation}\label{courb const}
    \beta:=- \frac{3}{4}(c_{12}^3)^2 - (c_{12}^1)^2 - (c_{12}^2)^2 + c_{12}^3 \, c_{23}^1 - X_2(c_{12}^1) +
    X_1(c_{12}^2) = const,
\end{equation}
by (\ref{spin metric}). We modify $(X_1,\,X_2,\,Z)$ by a rotation about $Z$. We obtain $(Y_1,\,Y_2,\,Z)$, where
$$ Y_1=\cos \psi \, X_1 + \sin \psi \, X_2 \quad \textrm{ and } \quad Y_2=-\sin \psi \, X_1 + \cos \psi \, X_2. $$
These vector fields satisfy
\begin{displaymath}
    \left\{ \begin{array}{rcl} \left[ Y_1,Y_2 \right] &=& (c_{12}^1 - X_1(\psi)) \, X_1 +
    (c_{12}^2 - X_2(\psi)) \, X_2 + c_{12}^3 \, Z \\ [2mm] \left[ Y_2,Z \right] &=& (Z(\psi) + c_{23}^1)
    Y_1 \\ [2mm] \left[ Z,Y_1 \right] &=& (Z(\psi) + c_{23}^1) Y_2.
    \end{array} \right.
\end{displaymath}

To construct a top, we need to find a function $\psi$ on $M$ such that $(Y_1,\,Y_2,\,Z)$ satisfies
(\ref{crochets}), that is, we have to integrate the following system~:
\begin{eqnarray*}
\left\{ \begin{array}{rcl}
    X_1(\psi) &=& c_{12}^1 \\
    X_2(\psi) &=& c_{12}^2 \\
    Z(\psi)&=& -c_{23}^1 + h,
\end{array} \right.
\end{eqnarray*}
where $h$ is a real number, such that $\beta=c_{12}^3\,h-\frac{3}{4}(c_{12}^3)^2$.

Let $\alpha$ be the differential 1-form defined by $\alpha(X_1)=c_{12}^1$, $\alpha(X_2)=c_{12}^2$ and
$\alpha(X_3)=-c_{23}^1+h$. Then the above system of PDEs can be written as
\begin{eqnarray} \label{equ D} d\psi=\alpha \end{eqnarray}
and there is a local solution if and only if $d\alpha=0$. This condition is equivalent to the first two equations
of (\ref{Jacobi2}) and to (\ref{courb const}). Thus, there is a local solution, hence a local top. Moreover, two
local solutions $\psi_1$ and $\psi_2$ differ only by a constant rotation, since we have
$d(\psi_1-\psi_2)=\alpha-\alpha=0$. So this local top is unique up to a reparametrization of its defining family
of 1-forms.

If $h=0$, this local top is integrable; if $h \neq 0$, it is a contact top. If $c_{12}^3=0$, we can choose $h$
arbitrarily. In that case, we can locally construct both integrable tops and contact tops.

There is a global solution of equation (\ref{equ D}) if and only if these local tops can be connected in a unique
way, that is, if the integral of $\alpha$ over every closed cycle in $M$ is an integer multiple of $2\pi$. This is
true under our assumption that $H^1_{de Rham}(M) = 0$. \finthm

\section{Metrics compatible with a given contact top} \label{compatible metrics}
 Let us now see how much freedom we have to choose a metric for a given top. \emph{On a manifold $(M,g)$ with a
given top, are there other metrics for which the same indexed pencil bundle defines a top?}

We are going to answer this question for contact tops only, because they are defined by round contact circles.
This determines a privileged plane bundle transverse to the axis bundle of the top, generated by the Reeb vector
fields associated to the generating forms of the contact circle. The roundness property implies that this plane
does not depend on the choice of these forms. We have the following result~:

\begin{prop} --- \label{toupie-métr}
Let $\mathcal{T}$ be a contact top on a Riemannian $3$-manifold $(M,g)$. Then $\mathcal{T}$ is a contact top for
another metric $g'$ if and only if the transition matrix between an orthonormal frame for $g$ and an orthonormal
frame for $g'$ is of the following form~:
$$  \left( \begin{array}{ccc}
    \lambda\,\alpha & \mu\,\gamma & 0 \\
    \lambda\,\beta & \mu\,\delta & 0 \\
    0 & 0 & \nu \\
    \end{array} \right) ,
\textrm{with}
    \left( \begin{array}{ccc}
    \alpha & \gamma \\
    \beta & \delta \\
    \end{array} \right)
\in {\rm O}(2,\R) \;\, \textrm{and} \;\, \lambda,\mu,\nu\in\R.$$
\end{prop}

\noindent {\bf Proof~:} Let $(X_1,\,X_2,\,X_3)$ be an orthonormal moving frame for $g$ which defines
$\mathcal{T}$ and thus satisfies (\ref{crochets}) for some constants $c$ and $k$. As $\mathcal{T}$ is a contact
top, we have $k\neq 0$. Then another moving frame $(Y_1,\,Y_2,\,Y_3)$ which defines the same indexed pencil
bundle, can be written as
$$ Y_1 = \lambda \, (\alpha \, X_1 + \beta \, X_2), \quad\;
Y_2 = \mu \, (\gamma \, X_1 + \delta \, X_2), \quad\; Y_3 = \nu \, X_3, $$ where $\alpha, \, \beta, \, \gamma$
and $\delta$ are real numbers such that \mbox{$\Delta:=\alpha \delta-\beta \gamma \neq 0$}, and where $\lambda$,
$\mu$ and $\nu$ are non-vanishing functions on $M$. Indeed, as $\mathcal{T}$ is defined by a taut, thus round,
contact circle (by Corollary \ref{cor str Cartan}), $Y_1$ and $Y_2$ span the same plane as $X_1$ and $X_2$.

As $(X_1,\,X_2,\,X_3)$ satisfies (\ref{crochets}), we have
\begin{eqnarray*}
    \left[ Y_1,Y_2\right]  &=& -\frac{Y_2(\lambda)}{\lambda}\,Y_1+\frac{Y_1(\mu)}{\mu}\,Y_2+c \,\Delta \frac{\lambda \,
    \mu}{\nu} \, Y_3 \\ [2mm]
    \left[ Y_2,Y_3\right] &=& \frac{\mu\,\nu\,(\gamma^2+\delta^2)}{\lambda\,\Delta}\,k\,Y_1- (\frac{Y_3(\mu)}{\mu}+
    \frac{\alpha\,\gamma+\beta\,\delta}{\Delta}\,\nu\,k)Y_2+\frac{Y_2(\nu)}{\nu}\,Y_3 \\[2mm]
    \left[ Y_3,Y_1\right] &=& (\frac{Y_3(\lambda)}{\lambda}-\frac{\alpha\,\gamma+\beta\,\delta}{\Delta}\,\nu\,k)Y_1+
    \frac{\lambda\,\nu\,(\alpha^2+\beta^2)}{\mu\,\Delta}\,k\,Y_2-\frac{Y_1(\nu)}{\nu}\,Y_3.
\end{eqnarray*}

Thus, $(Y_1,\,Y_2,\,Y_3)$ satisfies (\ref{crochets}) for some constants $\tilde{c}$ and $\tilde{k}$ if and only
if we have
$$ \begin{array}{l}
\displaystyle Y_1(\mu)=Y_2(\lambda)=Y_1(\nu)=Y_2(\nu)=0, \quad c\frac{\lambda\,\mu}{\nu}=const, \quad \frac{\mu}{\lambda}
\,(\gamma^2+\delta^2)=\frac{\lambda}{\mu}\,(\alpha^2+\beta^2), \\[1mm]
\displaystyle \frac{\mu\,\nu}{\lambda}=const, \quad \frac{\lambda\,\nu}{\mu}=const, \quad \mu\,\nu\,k(\alpha\,\gamma+
\beta\,\delta)+\Delta Y_3(\mu)=0, \\
\displaystyle \lambda\,\nu\,k(\alpha\,\gamma+\beta\,\delta)-\Delta Y_3(\lambda)=0.
\end{array} $$

So $\lambda$, $\mu$ and $\nu$ are constant on $M$. Up to multiplying $\lambda$ and $\mu$ by constant factors, we
can assume that $\alpha^2+\beta^2=\gamma^2+\delta^2=1$. The above relations show that the matrix $\left(
\begin{array}{cc} \alpha \; \gamma \\ \beta \; \delta \end{array} \right)$ is an orthogonal matrix. Thus,
$(Y_1,\,Y_2,\,Y_3)$ satisfies (\ref{crochets}) if and only if the transition matrix of $(Y_1,\,Y_2,\,Y_3)$ with
respect to $(X_1,\,X_2,\,X_3)$ corresponds to an endomorphism of the announced type. If $g'$ is the metric for
which $(Y_1,Y_2,Y_3)$ is orthonormal, then $\mathcal{T}$ is a top for $g'$, by Theorem \ref{car toupies}. \finthm

\section{Classification problems} \label{classification} 
Another natural problem is the classification of tops on a given Riemannian manifold. We will consider the
question of uniqueness~:

{\it Can two orthonormal global moving frames on a Riemannian 3-manifold $(M,g)$ determine two different tops,
that is, tops which are not related by an isometry or a reparametrization of the defining family of 1-forms~?}

\begin{lemma} \label{repar} ---
Two tops which have the same axis bundle and the same spinning direction are related by a reparametrization of
the defining family of 1-forms.
\end{lemma}

{\bf Proof :} Let $\mathcal{T}_1$ and $\mathcal{T}_2$ be two tops determined by two positively oriented
orthonormal moving frames $(X_1,X_2,Z)$ and $(Y_1,Y_2,Z)$, respectively, and let the associated constants be
$c,k$ and $\tilde{c},\tilde{k}$. By (\ref{courb extr}), the metric determines these constants up to simultaneous
multiplication by $-1$, so we have $c=\varepsilon \tilde{c}$ and $k=\varepsilon\tilde{k}$, with $\varepsilon=\pm
1$. The rotation speed of $\mathcal{T}_1$ and of $\mathcal{T}_2$ about their axis bundle along a geodesic
$\gamma$ which forms an angle $\varphi$ with the common axis bundle is given by
$$R_\gamma(X_1,Z)=(k-\frac{c}{2})\cos\varphi \quad \textrm{and} \quad
R_\gamma(Y_1,Z)=\varepsilon\,(k-\frac{c}{2})\cos\varphi.$$ As the spinning direction is the same, it follows that
$\varepsilon=1$.

So the rotation speeds along all geodesics are equal, and the two tops can only differ by a fixed rotation around
the axis bundle. This means that the second factors of the corresponding trivializations $\tau_1,\tau_2~:
\mathcal{F} \longrightarrow M\times\mathbb{P}^1(\mathbb{R})$ differ by a constant translation in
$\mathbb{P}^1(\mathbb{R})$, which can be expressed by a reparametrization of the defining family of 1-forms.
\finthm

\noindent {\bf Remark} --- It follows from this proof that two moving frames which determine two tops with the
same axis bundle define the same orientation if and only if their associated constants $c$ and $k$ coincide.
\vspace{4mm}

\noindent We will consider different situations with respect to the above uniqueness question~:
\begin{itemize}
\item If $(M,g)$ is not a space of constant curvature, Lemma \ref{Z unique} implies the uniqueness of the axis bundle
of a top on $M$. In that case, Lemma \ref{repar} grants the uniqueness of a top on $M$, up to reparametrization and
spinning direction.

\item On $\sph^3$ with the usual metric, any pivot field defines a contact top, by Theorem \ref{spin thm} and the
construction done in its proof. Furthermore, as two tops which are given by the same pivot field have the same axis
bundle, they coincide up to reparametrization and spinning direction. So we have to determine whether there are geodesic
Killing vector fields on $\sph^3$ which are not isometric to each other.

Let $\sph^3$ be the unit sphere of the space of quaternions $\ha$. Its algebra of Killing fields is of dimension
6 and it is generated by the vector fields which to a point $q$ associate the elements
  $qi$, $qj$, $qk$, $iq$, $jq$ and $kq$
of $T_q \sph^3$, respectively.

Let $v=a_1 \, qi + a_2 \, qj + a_3 \, qk + b_1 \, iq + b_2 \, jq + b_3 \, kq $ be a Killing vector field. Its
integral curves are geodesics if and only if either $a_1=a_2=a_3=0$ or $b_1=b_2=b_3=0$. Indeed, as the first
three vector fields $qi$, $qj$ and $qk$ define a moving frame, any vector field can be written as
$\omega=a\,qi+b\,qj+c\,qk$ with functions $a,b$ and $c$ on $M$. $\omega$ is geodesic if and only if $a,b$ and $c$
are constants. The same argument holds for the last three vector fields.

Hence, a pivot field $Z$ of the usual metric on $\sph^3$ is given either by \mbox{$a_1 \, qi + a_2 \, qj + a_3 \,
qk$} or by $b_1 \, iq + b_2 \, jq + b_3 \, kq$, where the $a_i$ and $b_i$ are real numbers. Two unitary vector
fields of one of these families are of course isometric. To see whether two vector fields of different families
are isometric, it is enough to consider the question for the vector fields $X_q = qi$ and $Y_q = iq$. A map $f$
exchanges $X$ and $Y$ if $Tf_q(X_q) = Y_{f(q)}$ for any point $q$ of $\sph^3$. For
$f(q_1,q_2,q_3,q_4)=(q_1,q_2,-q_3,-q_4)$, the tangent map $Tf$ satisfies this condition.

We have proved the following uniqueness result~:
\end{itemize}
\begin{prop} ---
On $\sph^3$ with the usual metric, all tops are isometric, up to spinning direction and reparametrization of the
defining family of 1-forms.
\end{prop}

\begin{itemize}
\item On a space $M$ of constant zero curvature, there might be non-isometric tops. In
particular, an isometry transforms a closed curve into a closed curve. On $\tor^3$, there are contact tops with
closed integral curves and others whose integral curves are not closed.
\end{itemize}

\begin{example} --- (Contact top on $\tor^3$ with closed integral curves.) \label{ex closed curves}

On $\tor^3$ with the pseudo-coordinates $(\theta_1,\,\theta_2,\,\theta_3)$, we consider the contact pencil bundle
associated to the contact circle of Example \ref{ex tore}, for $n=1$, generated by
$$ \omega_1=\cos\theta_1 \, d\theta_2+\sin\theta_1 \, d\theta_3 \quad \textrm{and} \quad
\omega_2=-\sin\theta_1 \, d\theta_2+\cos\theta_1 \, d\theta_3. $$

An adapted moving frame is given by the Reeb vector fields of $\omega_1$ and $\omega_2$ and by a third vector
field which defines the axis bundle~:
$$  R_1 = \cos\theta_1\,\frac{\partial}{\partial \theta_2} + \sin\theta_1\,\frac{\partial}{\partial \theta_3}, \quad
    R_2 = -\sin\theta_1\,\frac{\partial}{\partial \theta_2} + \cos\theta_1\,\frac{\partial}{\partial \theta_3},
\quad     X = \frac{\partial}{\partial \theta_1}.$$ The Lie brackets are $   \left[ R_1,\, R_2 \right] = 0, \,
\left[ R_2,\, X \right] = R_1$ and $\left[ X,\, R_1 \right]=R_2,$ so this moving frame determines a contact top on
$\tor^3$. The integral curves of its axis bundle are the integral curves of $X$, that is, the curves defined by
$\theta_2=\theta_3=const$, which are closed.
\end{example}

\begin{example} --- (Contact top on $\tor^3$ with non-closed integral curves.) \label{ex non-closed curves}

We modify a little the above example to get a moving frame of $\tor^3$ which defines the same metric, but
determines a different top, whose integral curves are not closed. Let us consider the following moving frame~:
\begin{eqnarray*}
\begin{array}{rcl}
    R_1 &=& \cos(\theta_1+ \varepsilon \,\theta_2)\,(\frac{\partial}{\partial \theta_2} - \varepsilon\,
    \frac{\partial}{\partial \theta_1}) + \sin(\theta_1+ \varepsilon\, \theta_2)\,\frac{\partial}{\partial \theta_3} \\ [1mm]
  R_2 &=& -\sin(\theta_1+ \varepsilon\, \theta_2)\,(\frac{\partial}{\partial \theta_2} - \varepsilon\,
  \frac{\partial}{\partial \theta_1}) + \cos(\theta_1+ \varepsilon\, \theta_2)\,\frac{\partial}{\partial \theta_3} \\ [1mm]
  X &=& \frac{\partial}{\partial \theta_1} + \varepsilon \,\frac{\partial}{\partial \theta_2}.
  \end{array}
\end{eqnarray*}
If $\varepsilon$ is not a rational number, the integral curves of $X$ are not closed, but each one is dense on a
2-torus contained in $\tor^3$. The Lie brackets are the following~:
$$  \left[ R_1,\, R_2 \right] = 0, \quad
    \left[ R_2,\, X \right] = (1+\varepsilon^2) \, R_1, \quad
    \left[ X,\, R_1 \right] = (1+\varepsilon^2) \, R_2, $$
so this moving frame determines a contact top. It is defined by the round contact circle generated by the forms
\begin{eqnarray*}
    \omega_1&=&\cos(\theta_1+ \varepsilon \,\theta_2) \, d(\theta_2- \varepsilon \,\theta_1)+\sin(\theta_1+ \varepsilon
    \,\theta_2) \, d\theta_3 \quad \textrm{and}\\
  \omega_2&=&-\sin(\theta_1+ \varepsilon \,\theta_2) \, d(\theta_2- \varepsilon \,\theta_1)+\cos(\theta_1+ \varepsilon
  \,\theta_2) \, d\theta_3.
\end{eqnarray*}
\end{example}

\section{Tops and Sasakian geometry}
The study of contact tops can be placed in the context of metric contact geometry. Therefore, it is suitable to
examine the relations to Sasakian geometry.

Let $(X_1,X_2,X_3)$ be a moving frame which determines a top on $(M,g)$ and thus satisfies (\ref{crochets}). Let
$(\omega_1,\omega_2,\omega_3)$ be the dual 1-forms. Then for $c \neq 0$, $\omega_3$ defines a K-contact
structure, that is, the Reeb vector field $X_3$ is Killing. In dimension 3, a manifold is Sasakian if and only if
it is K-contact (\cite{TS}, see also \cite{Bl}, Cor. 6.3 and 6.5), so in this case $(M,g)$ is Sasakian. On the
other hand, according to Itoh (see \cite{It}), no torus can carry a K-contact structure. But we have examples of
contact tops on $\tor^3$ (see Examples \ref{ex closed curves} and \ref{ex non-closed curves}), which correspond
to the case $c=0$.

So some tops define Sasakian structures and some do not, and the other way round a Sasakian structure $(\Phi,
\xi, \eta, g)$ (see \cite{Bl} for the notations) on a $3$-manifold $M$ defines a top in some cases and in others
does not. Indeed, according to R. Lutz (\cite{L}, §1.8), we can choose two vector fields $X_1$ and $X_2$ on $M$
such that $(X_1,\,X_2,\,\xi)$ is an orthonormal moving frame satisfying the Lie bracket relations~:
$$ \left[ X_1,X_2\right]=c_{12}^1 X_1+c_{12}^2 X_2+2 \xi, \quad \left[ X_2,\xi\right]=c_{23}^1 X_1,
\quad \left[ \xi,X_1\right]=c_{23}^1 X_2, $$ for some functions $c_{12}^1, \, c_{12}^2$ and $c_{23}^1$.

In view of expression (\ref{spin metric}), $g$ is then a spinning metric if and only if $$  - (c_{12}^1)^2 -
(c_{12}^2)^2 + 2 \, c_{23}^1 - X_2(c_{12}^1) + X_1(c_{12}^2) $$ is constant on $M$, that is, if the Sasakian
structure has constant sectional curvature. This is a necessary condition for the existence of a top on $M$,
according to Theorem~\ref{spin thm}. Furthermore, the Lie bracket relations imply that $\xi$ is a pivot field for
this metric. Thus, locally there are vector fields $\tilde{X_1}$ and $\tilde{X_2}$, such that
$(\tilde{X_1},\,\tilde{X_2},\,\xi)$ determines a top, as shown in the proof of Theorem \ref{spin thm}.
\vspace{4mm}

{\small Acknowledgement: This paper is based on my PhD-Thesis, Mulhouse 2004. I am very grateful to Robert Lutz
for our numerous discussions about tops and many more subjects. I also wish to thank Hansjörg Geiges for several
enlightening discussions. Furthermore, I am indebted to Hans Duistermaat for his constructive criticisms and many
important suggestions. During this work, I was partially supported by the SNF-project NMA1501 and by the grant GE
1245/1-2 within the DFG-Schwerpunktprogramm 1154 "Globale Differentialgeometrie".}

\end{document}